\documentclass[a4paper,10pt]{article}
\usepackage{stmaryrd}
\usepackage{amsfonts}
\usepackage{bbm}
\usepackage{amscd}
\usepackage{mathrsfs}
\usepackage{latexsym,amssymb,amsmath,amscd,amscd,amsthm,amsxtra,xypic}
\usepackage[dvips]{graphicx}
\usepackage[utf8]{inputenc}
\usepackage[T1]{fontenc}
\usepackage{lmodern}
\usepackage{amssymb}
\usepackage[all]{xy}
\usepackage{nicefrac,mathtools,enumitem}
\usepackage{microtype}
\usepackage{lmodern}

\newcommand {\emptycomment}[1]{}
\newtheorem{thm}{Theorem}[section]

\newtheorem{pro}[thm]{Proposition}

\newtheorem{rmk}[thm]{Remark}
\newtheorem{defi}[thm]{Definition}

\newcommand{\lon }{\,\rightarrow\,}
\newcommand{\be }{\begin{equation}}
\newcommand{\ee }{\end{equation}}

\newcommand{\pf}{\noindent{\bf Proof.}\ }

\newcommand{\g}{\frkg}

\newcommand{\Real}{\mathbb R}

\newcommand{\CWM}{C^{\infty}(M)}

\newcommand{\frkg}{\mathfrak g}
\newcommand{\frkh}{\mathfrak h}

\newcommand{\frkk}{\mathfrak k}
\newcommand{\frkl}{\mathfrak l}

\def\qed{\hfill ~\vrule height6pt width6pt depth0pt}

\newcommand{\br}[1]{   [ \cdot,    \cdot  ]   }

\newcommand{\gl}{\mathfrak {gl}}

\begin{document}
\title{
{ Equivalent description of Hom-Lie algebroids } }
\author{ Zhen Xiong \\
Department of Mathematics and Computer, Yichun University,
 Jiangxi 336000,  China
}
\date{email:205137@jxycu.edu.cn}
\footnotetext{{\it{Keyword}: Hom-Lie algebroids; Hom-Lie algebras; representations
 }}
 \footnotetext{{\it{MSC}}: 17B99,58H05}
 \footnotetext{Supported by the NSF of China (No.11771382) and the Science and Technology Project(GJJ161029)of Department of Education, Jiangxi Province.}

\maketitle

\begin{abstract}
 In this paper, we study representations of Hom-Lie algebroids, give some properties of  Hom-Lie algebroids and discuss equivalent statements of Hom-Lie algebroids. Then, we prove that two known definitions of Hom-Lie algebroids can be transformed into each other under some conditions.
\end{abstract}


\section{Introduction}

The notion of Hom-Lie algebras was introduced by Hartwig, Larsson,
and Silvestrov in \cite{HLS} as a part of a study of deformations of
the Witt and the Virasoro algebras. In a Hom-Lie algebra, the Jacobi
identity is twisted by a linear map, called the Hom-Jacobi identity.
Some $q$-deformations of the Witt and the Virasoro algebras have the
structure of a Hom-Lie algebra \cite{HLS,hu}. Because of close relations
to discrete and deformed vector fields and differential calculus
\cite{HLS,LD1,LD2},   more people pay special attention to this algebraic structure.
For a party of $k$-cochains on Hom-Lie algebras, name $k$-Hom-cochains, there is a series of coboundary operators \cite{yx}; for
regular Hom-Lie algebras, \cite{cy} gives a new coboundary operator on $k$-cochains, and there are many works have been done by the special coboundary operator \cite{cy,cl}. In \cite{xz}, there is a series of coboundary operators, and the author
 generalizes  the
result
" If $\frkk$ is a Lie algebra, $\rho:\frkk\longrightarrow\gl(V)$ is
a representation if and only if there is a degree-$1$ operator $D$
on $\Lambda\frkk^*\otimes V$ satisfying $D^2=0$, and
  $$
  D(\xi\wedge \eta\otimes u)=d_\frkk\xi\wedge\eta\otimes u+(-1)^k\xi \wedge D(\eta\otimes u),\quad \forall~\xi\in \wedge^k\frkk^*,~ \eta\in\wedge^l\frkk^*,~u\in V,
  $$
  where $d_\frkk:\wedge^k\g^*\longrightarrow \wedge^{k+1}\g^*$ is the coboundary operator associated to the trivial
  representation."

Geometric generalizations of Hom-Lie algebras are given in \cite{hom-Lie algebroids}\cite{cl}. In \cite{hom-Lie algebroids}, C. Laurent-Gengoux and J. Teles proved that there is a one-to-one correspondence between Hom-Gerstenhaber algebras and Hom-Lie algebroids; in \cite{QHC}, base on Hom-Lie algebroids from \cite{hom-Lie algebroids}, the authors study representation of Hom-Lie algebroids. In \cite{cl}, the authors make small modifications to the definition of Hom-Lie algebroids, and give a new definition of Hom-Lie algebroids, base on the new definition of Hom-Lie algebroids, definitions of Hom-Lie bialgebroids and Hom-Courant algebroids are given.

In this article, we first study representations of Hom-Lie algebroids, give equivalent statements of Hom-Lie algebroids and prove that different definitions of Hom-Lie algebroids are given by the same Hom-Lie algebras and their representations.

The paper is organized as follows. In Section 2, we recall some
basic notions. In Section 3, first, we study representations of Hom-Lie algebroids, and give some properties of Hom-Lie algebroids. Then, we prove that two known definitions of Hom-Lie algebroids can be transformed into each other(Theorem\ref{thm4}, Theorem\ref{thm5}).

\section{Preliminaries}
\subsection{Hom-Lie algebras and their representations}
The notion of a Hom-Lie algebra was introduced in \cite{HLS}, see also \cite{BM,MS2} for more information.
\begin{defi}
\begin{itemize}
\item[\rm(1)]
  A Hom-Lie algebra is a triple $(\frkg,\br ,,\alpha)$ consisting of a
  vector space $\frkg$, a skewsymmetric bilinear map (bracket) $\br,:\wedge^2\frkg\longrightarrow
  \frkg$ and a linear transformation $\alpha:\frkg\lon\frkg$ satisfying $\alpha[x,y]=[\alpha(x),\alpha(y)]$, and the following Hom-Jacobi
  identity:
  \begin{equation}
   [\alpha(x),[y,z]]+[\alpha(y),[z,x]]+[\alpha(z),[x,y]]=0,\quad\forall
x,y,z\in\frkg.
  \end{equation}

 A Hom-Lie algebra is called a regular Hom-Lie algebra if $\alpha$ is
a linear automorphism.

 \item[\rm(2)] A subspace $\frkh\subset\frkg$ is a Hom-Lie sub-algebra of $(\frkg,\br ,,\alpha)$ if
 $\alpha(\frkh)\subset\frkh$ and
  $\frkh$ is closed under the bracket operation $\br,$, i.e. for all $ x,y\in\frkh$,
  $[x,y] \in\frkh.  $
  \item[\rm(3)] A morphism from the  Hom-Lie algebra
$(\frkg,[\cdot,\cdot]_{\frkg},\alpha)$ to the hom-Lie algebra
$(\frkh,[\cdot,\cdot]_{\frkh},\gamma)$ is a linear map
$\psi:\frkg\longrightarrow\frkh$ such that
$\psi([x,y]_{\frkg})=[\psi(x),\psi(y)]_{\frkh}$ and
$\psi\circ \alpha =\gamma\circ \psi$.
  \end{itemize}
\end{defi}

 Representation and cohomology theories of Hom-Lie algebra are
systematically introduced in \cite{AEM,homlie1}. See \cite{Yao2} for
homology theories of Hom-Lie algebras.

\begin{defi}
  A representation of the Hom-Lie algebra $(\frkg,\br,,\alpha)$ on
  a vector space $V$ with respect to $\beta\in\gl(V)$ is a linear map
  $\rho:\frkg\longrightarrow \gl(V)$, such that for all
  $x,y\in\frkg$, the following equalities are satisfied:
  \begin{eqnarray*}
 \label{representation1} \rho (\alpha(x))\circ \beta&=&\beta\circ \rho (x);\\
    \rho([x,y] )\circ
    \beta&=&\rho (\alpha(x))\circ\rho (y)-\rho (\alpha(y))\circ\rho (x).
  \end{eqnarray*}
\end{defi}
Let $(\frkg,[\cdot,\cdot],\alpha)$ be a Hom-Lie algebra, $V$ be a
vector space, $\rho:\frkg\longrightarrow\frkg\frkl(V)$ be a
representation of $(\frkg,[\cdot,\cdot],\alpha)$ on the vector space
$V$ with respect to $\beta\in \frkg\frkl(V)$.

The set of {\bf
$k$-cochains} on $\frkg$ with values in $V$, which we denote by
$C^k(\frkg;V)$, is the set of skewsymmetric $k$-linear maps from
$\frkg\times\cdots\times\frkg$($k$-times) to $V$:
$$C^k(\frkg;V):=\{\eta:\wedge^k\frkg\longrightarrow V ~\mbox{is a
linear map}\}.$$

In \cite{xz}, when $\beta\in GL(V)$, there is a series operators $d^s:C^k(\frkg;V)\longrightarrow
C^{k+1}(\frkg;V)$ is given by
\begin{eqnarray*}
d^s\eta(x_1,\cdots,x_{k+1})&=&\sum_{i=1}^{k+1}(-1)^{i+1}\beta^{k+1+s}\rho(x_i)\beta^{-k-2-s}\eta(\alpha(x_1),\cdots,\hat{x_i},\cdots,\alpha(x_{k+1}))\\
&&+\sum_{i<j}(-1)^{i+j}\eta([x_i,x_j],\alpha(x_1),\cdots,\widehat{x_{i,j}},\cdots,\alpha(x_{k+1})),
\end{eqnarray*}
where $\beta^{-1}$ is the inverse of $\beta$, $\eta\in C^k(\frkg;V)$, and the author have the results: $d^s\circ d^s=0$.

\subsection{Hom-Lie algebroids}
Now, we introduce two kinds of definitions of Hom-Lie algebroids, they are from \cite{hom-Lie algebroids} and \cite{cl}. More about Hom-Lie algebroids, please see \cite{hom-Lie algebroids} and \cite{cl}.
\begin{defi}\label{s1-1}\cite{hom-Lie algebroids} A Hom-Lie algebroid is a quintuple
$(A,\varphi,[\cdot,\cdot],\rho_A,\alpha_A)$, where $A$ is a vector
bundle over a manifold $M, \varphi:M\longrightarrow M$ is a smooth
map, $[\cdot,\cdot]:\Gamma( A)\otimes \Gamma( A) \longrightarrow
 \Gamma( A)$ is a bilinear map,called bracket, $\rho_A:\varphi^!A\longrightarrow\varphi^!TM$ is a vector
 bundle morphism, called anchor and $\alpha_A:\Gamma( A) \longrightarrow
 \Gamma( A)$ is a linear endomorphism of $\Gamma( A)$, for $
X,Y\in\Gamma( A),f\in\CWM$ such that:
\begin{itemize}
\item[(1)] $\alpha_A(fX)=\varphi^*(f)\alpha_A(X)$;
\item[(2)] the triple $(\Gamma( A),[\cdot,\cdot],\alpha_A)$ is a Hom-Lie
algebra;
\item[(3)] the following Hom-Leibniz identity holds:
$$[X,fY]=\varphi^*(f)[X,Y]+\rho_A(X)(f)\alpha_A(Y);$$
\item[(4)] $\rho_A$ is a representation of Hom-Lie algebra $(\Gamma(
A),[\cdot,\cdot],\alpha_A)$ on $\CWM$ with respect to $\varphi^*$.
\end{itemize}
\end{defi}
In fact, according to Definition\ref{s1-1}, for $ X,Y\in\Gamma(
A),f,g\in\CWM$, we have the following properties:
\begin{itemize}
\item[(a)] $\alpha_A([X,Y])=[\alpha_A(X),\alpha_A(Y)];$
\item[(b)] $\varphi^*:\CWM\longrightarrow\CWM
$, defined by $ \varphi^*(f)=f\circ\varphi,\varphi^*$ is a morphism
of $ \CWM $ ;
\item[(c)] $\rho_A(X)(fg)=\varphi^*(f)\rho_A(X)(g)+\varphi^*(g)\rho_A(X)(f)$.According to (3) in Definition:
$$\rho_A(X)(fg)=[[X,fg]]=[[X,f]]\alpha_A(g)+[[X,g]]\alpha_A(f)=\alpha_A(g)\rho_A(X)(f)+\alpha_A(f)\rho_A(X)(g)$$
where $[[\cdot ,\cdot]]$ is define in Definition 3.1 of \cite{hom-Lie algebroids}.
\item[(d)] $\alpha_A(f)=\varphi^*(f)$, when $\alpha=\rm{id}$, then
$\varphi^*=id$, Hom-Lie algebroid
$(A,\varphi,[\cdot,\cdot],\rho_A,\alpha_A)$ is just a Lie algebroid;
\item[(e)] $\rho_A(fX)=\varphi^*(f)\rho_A(X)$. It follows from:
$$\rho_A(fX)(g)=[[fX,g]]=\alpha_A(f)[[X,g]]=\varphi^*(f)\rho_A(X)(g).$$
\end{itemize}
\begin{defi}\label{s1-2}\cite{cl} A Hom-Lie algebroid is a quintuple
$(B,\varphi,[\cdot,\cdot],\rho_B,\alpha_B)$, where $B$ is a vector
bundle over a manifold $M, \varphi:M\longrightarrow M$ is a smooth
map, $[\cdot,\cdot]:\Gamma( B)\otimes \Gamma( B) \longrightarrow
 \Gamma( B)$ is a bilinear map,called bracket, $\rho_B:B\longrightarrow\varphi^!TM$ is
 a bundle map, called anchor and $\alpha_B:\Gamma( B) \longrightarrow
 \Gamma( B)$ is a linear endomorphism of $\Gamma( B)$, for $
X,Y\in\Gamma(B),f\in\CWM$ such that:
\begin{itemize}
\item[1)] $\alpha_B(fX)=\varphi^*(f)\alpha_B(X)$;
\item[2)] the triple $(\Gamma( B),[\cdot,\cdot],\alpha_B)$ is a Hom-Lie
algebra;
\item[3)] the following Hom-Leibniz identity holds:
$$[X,fY]=\varphi^*(f)[X,Y]+\rho_B(\alpha_B(X))(f)\alpha_B(Y);$$
\item[4)] $\rho_B$ is a representation of Hom-Lie algebra $(\Gamma(
B),[\cdot,\cdot],\alpha_B)$ on $\CWM$ with respect to $\varphi^*$.
\end{itemize}
\end{defi}
From Defintion \ref{s1-2}, for $ X,Y\in\Gamma(
A),f,g\in\CWM$, we have:
\begin{itemize}
\item[a)] $\alpha_B([X,Y])=[\alpha_B(X),\alpha_B(Y)];$
\item[b)] $\varphi^*:\CWM\longrightarrow\CWM
$, defined by $ \varphi^*(f)=f\circ\varphi,\varphi^*$ is a morphism
of $ \CWM $ ;
\item[c)] $\rho_B(X)(fg)=\varphi^*(f)\rho_B(X)(g)+\varphi^*(g)\rho_B(X)(f)$.
\item[d)] $\alpha_B(f)=\varphi^*(f)$.
\item[e)] $\rho_B(fX)=f\rho_B(X)$.
\end{itemize}
When $\alpha_B$ and $\varphi$ are invertible, Hom-Lie bialgebroids and Hom-Courant algebroids are given in \cite{cl}.

\section{Representations of Hom-Lie algebroids}
In this section, we assume that map $\varphi:M\longrightarrow M$ is a involution, i.e. $\varphi^2=\rm{id}$.

Let $(E,\varphi,[\cdot,\cdot],\rho_E,\alpha_E)$ be a Hom-Lie algebroid.
Whence $(\rho_E,\varphi^*)$ is a representation of $(\Gamma(
E),[\cdot,\cdot],\alpha_E)$ on $\CWM$. Where $E$ is $A$ or $B$.
We define $d^s:C^k(\Gamma(E);\CWM)\longrightarrow C^{k+1}(\Gamma(E);\CWM),\quad s=0,1,\ldots$, by setting
\begin{eqnarray*}
d^s\eta(X_1,\cdots,X_{k+1})&=&\sum_{i=1}^{k+1}(-1)^{i+1}{\varphi^*}^{k+1+s}\rho_E(X_i){\varphi^*}^{-k-2-s}\eta(\alpha_E(X_1),\cdots,\hat{X_i},\cdots,\alpha_E(X_{k+1}))\\
&&+\sum_{i<j}(-1)^{i+j}\eta([X_i,X_j],\alpha_E(X_1),\cdots,\hat{X_i},\cdots,\hat{X_j},\cdots,\alpha_E(X_{k+1})),
\end{eqnarray*}
where $C^k(\Gamma(E);\CWM)=\{\eta:\wedge^k\Gamma(E)\longrightarrow\CWM \quad \mbox{is a $\Real$ linear map}\}$, $X_i\in\Gamma(E).$

We define map $\alpha_E^*:C^k(\Gamma(E);\CWM)\longrightarrow C^k(\Gamma(E);\CWM)$ by
$$\alpha_E^*(\eta)(X_1,\cdots,X_k)=\varphi^*\circ\eta(\alpha_E(X_1),\cdots,\alpha_E(X_k)),\quad X_i\in\Gamma(E).$$
When $f\in C^0(\Gamma(E);\CWM)=\CWM$, we have $\alpha_E^*(f)=\varphi^*(f)$.

Let $\Gamma(\wedge^k E^*)=\{\eta\in C^k(\Gamma(E);\CWM) |\quad f\in\CWM,\quad \eta(fX_1,\cdots,X_k)=f\eta(X_1,\cdots,X_k)\}$.
Then, $\Gamma(\wedge^k E^*)$ is a subset of $C^k(\Gamma(E);\CWM)$. Let $\alpha_E^*$ acts on $\Gamma(\wedge^k E^*)$, we have:
$$\alpha_E^*:\Gamma(\wedge^k E^*)\longrightarrow \Gamma(\wedge^k E^*).$$
Actually, when $\xi\in\Gamma(E^*)$, for $f\in\CWM$, we have
$$
\alpha_E^*(\xi)(fX)=\varphi^*\circ\xi(\alpha_E(fX))=\varphi^*\circ\xi(\varphi^*(f)\alpha_E(X))={\varphi^*}^2(f)\varphi^*\circ\xi(\alpha_E(X))=f\alpha_E^*(\xi)(X),
$$
so, $\alpha_E^*(\xi)\in \Gamma(E^*)$, and we have
$$\alpha_E^*(\xi_1\wedge\cdots\xi_k)=\alpha_E^*(\xi_1)\wedge\cdots\alpha_E^*(\xi_k),\quad \xi_1,\cdots,\xi_k\in\Gamma(E^*).$$
Let $\Gamma(\wedge^\bullet E^*)=\oplus_k\Gamma(\wedge^k E^*)$, $C^\bullet(\Gamma(E);\CWM)=\oplus_kC^k(\Gamma(E);\CWM)$, we have:
$\Gamma(\wedge^\bullet E^*)$ is a subset of $C^\bullet(\Gamma(E);\CWM)$. So, if $\xi\in \Gamma(\wedge^\bullet E^*)$, we have $d^s\xi\in C^\bullet(\Gamma(E);\CWM)$.

At the same time, $\varphi^*$ can induce a map $\overline{\varphi^*}:C^k(\Gamma(E);\CWM)\longrightarrow C^k(\Gamma(E);\CWM)$, which define by
$$\overline{\varphi^*}(\eta)(X_1,\cdots,X_k)=\varphi^*\circ\eta(X_1,\cdots,X_k),\quad X_i\in\Gamma(E).$$
Then, we have:
$$\overline{\varphi^*}(\eta_1\wedge\cdots\eta_k)=\overline{\varphi^*}(\eta_1)\wedge\cdots\overline{\varphi^*}(\eta_k).$$

\begin{pro}
With above notations, for $\eta_1\in C^k(\Gamma(E);\CWM),\eta_2\in
C^l(\Gamma(E);\CWM)$, we have:
$$d^s(\eta_1\wedge\eta_2)=d^{s+l}\eta_1\wedge\overline{\varphi^*}\circ\alpha_E^*(\eta_2)+(-1)^k\overline{\varphi^*}\circ\alpha_E^*(\eta_1)\wedge
d^{s+k}\eta_2.$$
\end{pro}
\pf First let $k=1$, $\eta_1\wedge\eta_2\in C^{l+1}(\Gamma(E);\CWM)$, we have:
\begin{eqnarray*}
&&d^s(\eta_1\wedge\eta_2)(X_1,\cdots,X_{l+2})\\
&=&\sum_{i=1}(-1)^{i+1}{\varphi^*}^{l+2+s}\rho_E(X_i){\varphi^*}^{-l-3-s}\eta_1\wedge\eta_2(\alpha_E(x_1),\cdots,\hat{X_i},\cdots,\alpha_E(x_{l+2}))\\
&&+\sum_{i<j}(-1)^{i+j}\eta_1\wedge\eta_2([X_i,X_j],\alpha_E(X_1),\cdots,\widehat{X_{i,j}},\cdots,\alpha_E(X_{l+2}))\\
&=&\sum_{i<q}(-1)^{q+i+1}\eta_2(\alpha_E(X_1),\cdots,\widehat{X_{i,q}},\cdots,\alpha_E(X_{l+2})){\varphi^*}^{l+2+s}\rho_E(X_i){\varphi^*}^{-l-3-s}\eta_1(\alpha_E(X_q))\\
&&+\sum_{q<i}(-1)^{i+q}\eta_1(\alpha_E(X_q)){\varphi^*}^{l+2+s}\rho_E(X_i){\varphi^*}^{-l-3-s}\eta_2(\alpha_E(X_1),\cdots,\widehat{X_{q,i}},\cdots,\alpha_E(X_{l+2}))\\
&&+\sum_{i<j}(-1)^{i+j}\eta_1([X_i,X_j]\eta_2(\alpha_E(X_1),\cdots,\widehat{X_{i,j}},\cdots,\alpha_E(X_{l+2}))\\
&&+\sum_{i<q}(-1)^{q+i+1}\eta_1(\alpha_E(x_q)){\varphi^*}^{l+2+s}\rho_E(X_i){\varphi^*}^{-l-3-s}\eta_2(\alpha_E(X_1),\cdots,\widehat{X_{i,q}},\cdots,\alpha_E(X_{l+2}))\\
&&+\sum_{q<i}(-1)^{i+q}\eta_2(\alpha_E(X_1),\cdots,\widehat{X_{q,i}},\cdots,\alpha_E(X_{l+2})){\varphi^*}^{l+2+s}\rho_E(X_i){\varphi^*}^{-l-3-s}\eta_1(\alpha_E(X_q))\\
&&+\sum_{q<i<j}(-1)^{q+i+j}\eta_1(\alpha_E(X_q))\eta_2([X_i,X_j],\alpha_E(X_1),\cdots,\widehat{x_{q,i,j}},\cdots,\alpha_E(X_{l+2}))\\
&&+\sum_{i<q<j}(-1)^{q+i+j+1}\eta_1(\alpha_E(X_q))\eta_2([X_i,X_j],\alpha_E(X_1),\cdots,\widehat{x_{i,q,j}},\cdots,\alpha_E(X_{l+2}))\\
&&+\sum_{i<j<q}(-1)^{q+i+j}\eta_1(\alpha_E(X_q))\eta_2([X_i,X_j],\alpha_E(X_1),\cdots,\widehat{x_{i,j,q}},\cdots,\alpha_E(X_{l+2}))\\
&=&d^{s+l}\eta_1\wedge\overline{\varphi^*}\circ\alpha_E^*(\eta_2)(X_1,\cdots,X_{l+2})+(-1)\overline{\varphi^*}\circ\alpha_E^*(\eta_1)\wedge
d^{s+1}\eta_2(X_1,\cdots,X_{l+2}).
\end{eqnarray*}
So, when $k=1$, we have:
$$d^s(\eta_1\wedge\eta_2)=d^{s+l}\eta_1\wedge\overline{\varphi^*}\circ\alpha_E^*(\eta_2)+(-1)\overline{\varphi^*}\circ\alpha_E^*(\eta_1)\wedge
d^{s+1}\eta_2.$$
By induction on $k$,  assume that when $k=n$, we have:
$$d^s(\eta_1\wedge\eta_2)=d^{s+l}\eta_1\wedge\overline{\varphi^*}\circ\alpha_E^*(\eta_2)+(-1)^n\overline{\varphi^*}\circ\alpha_E^*(\eta_1)\wedge
d^{s+n}\eta_2.$$
For any $\eta_3\in C^1(\Gamma(E);\CWM)$, then $\eta_1\wedge\eta_3\in C^{n+1}(\Gamma(E);\CWM)$, we have
\begin{eqnarray*}
d^s(\eta_1\wedge\eta_3\wedge\eta_2)&=&d^s(\eta_1\wedge(\eta_3\wedge\eta_2))\\
&=&d^{s+l+1}\eta_1\wedge\overline{\varphi^*}\circ\alpha_E^*(\eta_3\wedge\eta_2)+(-1)^n\overline{\varphi^*}\circ\alpha_E^*(\eta_1)\wedge
d^{s+n}(\eta_3\wedge\eta_2)\\
&=&d^{s+l+1}\eta_1\wedge\overline{\varphi^*}\circ\alpha_E^*(\eta_3)\wedge\overline{\varphi^*}\circ\alpha_E^*(\eta_2)\\
&&+(-1)^n\overline{\varphi^*}\circ\alpha_E^*(\eta_1)\wedge\Big(d^{s+n+l}\eta_3\wedge\overline{\varphi^*}\circ\alpha_E^*(\eta_2)+(-1)\overline{\varphi^*}\circ\alpha_E^*(\eta_3)\wedge d^{s+n+1}\eta_2\Big)   \\
&=&\Big(d^{s+l+1}\eta_1\wedge\overline{\varphi^*}\circ\alpha_E^*(\eta_3)+(-1)^n\overline{\varphi^*}\circ\alpha_E^*(\eta_1)\wedge d^{s+n+l}\eta_3\Big)\wedge\overline{\varphi^*}\circ\alpha_E^*(\eta_2)\\
&&+(-1)^{n+1}\overline{\varphi^*}\circ\alpha_E^*(\eta_1)\wedge\overline{\varphi^*}\circ\alpha_E^*(\eta_3)\wedge d^{s+n+1}\eta_2\\
&=&d^{s+l}(\eta_1\wedge\eta_3)\wedge\overline{\varphi^*}\circ\alpha_E^*(\eta_2)+(-1)^{n+1}\overline{\varphi^*}\circ\alpha_E^*(\eta_3\wedge\eta_1)\wedge
d^{s+n+1}\eta_2.
\end{eqnarray*}
The proof is completed.\qed

\begin{pro}
With above notations, we have:
$$\alpha_E^*\circ d^s=d^s\circ\alpha_E^*;\quad \overline{\varphi^*}\circ d^s=d^{s+1}\circ\overline{\varphi^*}.$$
\end{pro}
\pf With straightforward computations, for any $\eta\in C^k(\Gamma(E);\CWM)$, we have:
\begin{eqnarray*}
&&\alpha_E^*\circ d^s\eta(X_1,\cdots,X_{k+1})\\
&=&\varphi^*\circ d^s\eta(\alpha_E(X_1),\cdots,\alpha_E(X_{k+1}))\\
&=&\sum_{i=1}(-1)^{i+1}{\varphi^*}^{k+2+s}\rho_E(\alpha_E(X_i)){\varphi^*}^{-k-2-s}\eta(\alpha_E^2(X_1),\cdots,\hat{X_i},\cdots,\alpha_E^2(X_{k+1}))\\
&&+\sum_{i<j}(-1)^{i+j}\varphi^*\circ\eta([\alpha_E(X_i),\alpha_E(X_j)],\alpha_E^2(X_1),\cdots,\widehat{X_{i,j}},\cdots,\alpha_E^2(X_{k+1}))\\
&=&\sum_{i=1}(-1)^{i+1}{\varphi^*}^{k+3+s}\rho_E(X_i){\varphi^*}^{-k-4-s}\alpha_E^*(\eta)(\alpha_E(X_1),\cdots,\hat{X_i},\cdots,\alpha_E(X_{k+1}))\\
&&+\sum_{i<j}(-1)^{i+j}\alpha_E^*(\eta)([X_i,X_j],\alpha_E(X_1),\cdots,\widehat{X_{i,j}},\cdots,\alpha_E(X_{k+1}))\\
&=&d^s\circ\alpha_E^*(\eta)(X_1,\cdots,X_{k+1}).
\end{eqnarray*}
At the same time, we have:
\begin{eqnarray*}
&&\overline{\varphi^*}\circ d^s\eta(X_1,\cdots,X_{k+1})\\
&=&\sum_{i=1}(-1)^{i+1}{\varphi^*}^{k+2+s}\rho_E(X_i){\varphi^*}^{-k-2-s}\eta(\alpha_E(X_1),\cdots,\hat{X_i},\cdots,\alpha_E(X_{k+1}))\\
&&+\sum_{i<j}(-1)^{i+j}\varphi^*\circ\eta([X_i,X_j],\alpha_E(X_1),\cdots,\widehat{X_{i,j}},\cdots,\alpha_E(X_{k+1}))\\
&=&\sum_{i=1}(-1)^{i+1}{\varphi^*}^{k+2+s}\rho_E(X_i){\varphi^*}^{-k-3-s}\overline{\varphi^*}(\eta)(\alpha_E(X_1),\cdots,\hat{X_i},\cdots,\alpha_E(X_{k+1}))\\
&&+\sum_{i<j}(-1)^{i+j}\overline{\varphi^*}(\eta)([X_i,X_j],\alpha_E(X_1),\cdots,\widehat{X_{i,j}},\cdots,\alpha_E(X_{k+1}))\\
&=&d^{s+1}\circ\overline{\varphi^*}(\eta)(X_1,\cdots,X_{k+1}).
\end{eqnarray*}
We complete this proof.\qed

Now, we revisited representations of Hom-Lie algebroids respectively base on Definition\ref{s1-1} and Definition\ref{s1-2}.

\begin{thm}\label{thm4}
Let  $A$ be a vector
bundle over  manifold $M, \varphi:M\longrightarrow M$ is a smooth
map and $\varphi^2=\rm{id}$, $\alpha_A:\Gamma( A) \longrightarrow
 \Gamma( A)$ is a linear endomorphism of $\Gamma( A)$ i.e.for $f\in\CWM,X\in\Gamma( A), \alpha_A(fX)=\varphi^*(f)\alpha_A(X)$. Then $(A,\varphi,[\cdot,\cdot],\rho_A,\alpha_A)$ is a Hom-Lie algebriod define by Definition\ref{s1-1} if and only if there is a
series operators $d^s:C^k(\Gamma(A);\CWM)\longrightarrow C^{k+1}(\Gamma(A);\CWM), s=0,1,\ldots$, and such that:
\begin{itemize}
\item[\rm{(i)}]$d^s\circ d^s=0$;
\item[\rm{(ii)}]for any $\eta_1\in C^k(\Gamma(A);\CWM),\eta_2\in C^l(\Gamma(A);\CWM)$, we have
$$d^s(\eta_1\wedge\eta_2)=d^{s+l}\eta_1\wedge\overline{\varphi^*}\circ\alpha_A^*(\eta_2)+(-1)^k\overline{\varphi^*}\circ\alpha_A^*(\eta_1)\wedge
d^{s+k}\eta_2.$$
\item[\rm{(iii)}]$\alpha_A^*\circ d^s=d^s\circ\alpha_A^*;\quad \overline{\varphi^*}\circ d^s=d^{s+1}\circ\overline{\varphi^*}.$
\item[\rm{(iv)}]for $f\in\CWM=C^0(\Gamma(A);\CWM)$, we have: $d^0f\in \Gamma(A^*)$.
\item[\rm{(v)}]for $f\in\CWM, \xi\in\Gamma(A^*)$, we have;
$$d^0\xi(fX_1,X_2)=\varphi^*(f)d^0\xi(X_1,X_2).$$
\end{itemize}
\end{thm}
\pf For necessity, with above Propositions which we proved, we just need to prove \rm{(iv)} and \rm{(v)}.
For Hom-Lie algebriod $(A,\varphi,[\cdot,\cdot],\rho_A,\alpha_A)$ and $f, g\in\CWM=C^0(\Gamma(E);\CWM), X\in\Gamma(A)$, by the definition of $d^s$, we have:
$$d^0f(gX)=\varphi^*\rho_A(gX)f=g\varphi^*\rho_A(X)f=gd^0f(X).$$

For $f\in\CWM, \xi\in\Gamma(A^*)$, we have:
\begin{eqnarray*}
d^0\xi(fX_1,X_2)&=&\rho_A(fX_1)\varphi^*\xi(\alpha_A(X_2))-\rho_A(X_2)\varphi^*\xi(\alpha_A(fX_1))\\
&&-\xi([fX_1,X_2])\\
&=&\varphi^*(f)\rho_A(X_1)\varphi^*\xi(\alpha_A(X_2))-\varphi^*(f)\rho_A(X_2)\varphi^*\xi(\alpha_A(X_1))-\xi(\alpha_A(X_1))\rho_A(X_2)f\\
&&-\xi(\varphi^*(f)[X_1,X_2]-\rho_A(X_2)f\alpha_A(X_1))\\
&=&\varphi^*(f)d^0\xi(X_1,X_2).
\end{eqnarray*}
So, we proved the necessity of this Theorem. Now, we prove the adequacy of this theorem.

Sept1, we define $\rho_A:\varphi^!A\longrightarrow\varphi^!TM$ by:
\begin{equation}\label{s30}
\varphi^*\rho_A(X)f=d^0f(X),\quad X\in\Gamma(A),f\in\CWM.
\end{equation}
Then, by \rm{(iv)}, for $g\in\CWM$, $\varphi^*\rho_A(gX)f=d^0f(gX)=gd^0f(X)=g\varphi^*\rho_A(X)f$, we have:
\begin{equation}\label{s31}
\rho_A(gX)=\varphi^*(g)\rho_A(X).
\end{equation}
On the other hand, we have:
$$d^0(fg)(X)=\varphi^*\rho_A(X)(fg).$$
By \rm{(ii)}, for $f,g\in C^0(\Gamma(A);\CWM)=\CWM$, we have:
$$d^0(fg)=d^0fg+fd^0g.$$
Then, we have:
\begin{equation}\label{s32}
\rho_A(X)(fg)=\varphi^*(g)\rho_A(X)f+\varphi^*(f)\rho_A(X)g.
\end{equation}
The definition of $\rho_A$ is reasonable.

By \rm{(iii)}, we have:
\begin{eqnarray*}
\alpha_A^*\circ d^0f(X)&=&\varphi^*\circ d^0f(\alpha_A(X))\\
&=&\varphi^*\circ\varphi^*\rho_A(\alpha_A(X))f\\
&=&d^0\alpha_A^*(f)(X)\\
&=&\varphi^*\circ\rho_A(X)\varphi^*(f).
\end{eqnarray*}
We find the result:
\begin{equation}\label{s33}
\rho_A(\alpha_A(X))\varphi^*=\varphi^*\rho_A(X).
\end{equation}
Sept2, for any $\xi\in\Gamma(A^*),X,Y\in\Gamma(A)$, we define $[\cdot,\cdot]:\Gamma(A)\wedge\Gamma(A)\longrightarrow\Gamma(A)$by
\begin{equation}\label{s34}
\xi([X,Y])=\rho_A(X)\varphi^*\xi(\alpha_A(Y))-\rho_A(Y)\varphi^*\xi(\alpha_A(X))-d^0\xi(X,Y).
\end{equation}
So, by (\ref{s34}), we have:
\begin{eqnarray*}
\alpha_A^*(\xi)([X,Y])&=&\rho_A(X)\varphi^*\alpha_A^*(\xi)(\alpha_A(Y))-\rho_A(Y)\varphi^*\alpha_A^*(\xi)(\alpha_A(X))-d^0\alpha_A^*(\xi)(X,Y)\\
&=&\rho_A(X)\xi(\alpha_A^2(Y))-\rho_A(Y)\xi(\alpha_A^2(X))-d^0\alpha_A^*(\xi)(X,Y).
\end{eqnarray*}
By $\alpha_A^*(\xi)([X,Y])=\varphi^*\xi(\alpha_A[X,Y])$, (\ref{s33})and \rm{(iii)}, (\ref{s34}), we have
\begin{eqnarray*}
\xi(\alpha_A([X,Y]))&=&\varphi^*\rho_A(X)\xi(\alpha_A^2(Y))-\varphi^*\rho_A(Y)\xi(\alpha_A^2(X))-\varphi^*\circ\alpha_A^*\circ d^0\xi(X,Y)\\
&=&\rho_A(\alpha_A(X))\varphi^*\xi(\alpha_A^2(Y))-\rho_A(\alpha_A(Y))\varphi^*\xi(\alpha_A^2(X))-d^0\xi(\alpha_A(X),\alpha_A(Y))\\
&=&\xi([\alpha_A(X),\alpha_A(Y)]).
\end{eqnarray*}
We have the following:
\begin{equation}\label{s35}
\alpha_A([X,Y])=[\alpha_A(X),\alpha_A(Y)].
\end{equation}
For any $f\in\CWM$, by \rm{(i)}, (\ref{s30}) and (\ref{s34}), we have:
\begin{eqnarray*}
0&=&d^0\circ d^0f(X,Y)\\
&=&\rho_A(X)\varphi^*d^0f(\alpha_A(Y))-\rho_A(Y)\varphi^*d^0f(\alpha_A(X))-d^0f([X,Y])\\
&=&\rho_A(X)\rho_A(\alpha_A(Y))f-\rho_A(Y)\rho_A(\alpha_A(X))f-\varphi^*\rho_A([X,Y])f.
\end{eqnarray*}
We get
\begin{eqnarray}
\rho_A(X)\rho_A(\alpha_A(Y))-\rho_A(Y)\rho_A(\alpha_A(X))&=&\varphi^*\rho_A([X,Y])\nonumber\\
\varphi^*\rho_A(X)\rho_A(\alpha_A(Y))\varphi^*-\varphi^*\rho_A(Y)\rho_A(\alpha_A(X))\varphi^*&=&\rho_A([X,Y])\varphi^*\nonumber\\
\rho_A(\alpha_A(X))\rho_A(Y)-\rho_A(\alpha_A(Y))\rho_A(X)&=&\rho_A([X,Y])\varphi^*.\label{s36}
\end{eqnarray}
Sept3, by (\ref{s34}), (\ref{s31}), (\ref{s32}) and \rm{(v)}, we have:
\begin{eqnarray*}
\xi([X,fY])&=&\rho_A(X)\varphi^*\xi(\alpha_A(fY))-\rho_A(fY)\varphi^*\xi(\alpha_A(X))-d^0\xi(X,fY)\\
&=&\rho_A(X)(f\varphi^*\xi(\alpha_A(Y)))-\varphi^*(f)\rho_A(Y)\varphi^*\xi(\alpha_A(X))-\varphi^*(f)d^0\xi(X,Y)\\
&=&\varphi^*(f)\xi([X,Y])+\xi(\rho_A(X)f\alpha_A(Y)).
\end{eqnarray*}
So, for $f\in\CWM, X,Y\in\Gamma(A)$, we have:
\begin{equation}\label{s37}
[X,fY]=\varphi^*(f)[X,Y]+\rho_A(X)f\alpha_A(Y).
\end{equation}
Sept4, by \rm{(iii)}, we have
$$d^1=\overline{\varphi^*}\circ d^0\circ\overline{\varphi^*}.$$
For $\eta_1,\eta_2\in C^1(\Gamma(A);\CWM)$, by \rm{(ii)}, we have
\begin{eqnarray*}
d^0(\eta_1\wedge\eta_2)&=&d^1\eta_1\wedge\overline{\varphi^*}\circ\alpha_A^*(\eta_2)-\overline{\varphi^*}\circ\alpha_A^*(\eta_1)\wedge d^1\eta_2\\
&=&\overline{\varphi^*}\Big(d^0\overline{\varphi^*}(\eta_1)\wedge\alpha_A^*(\eta_2)-\alpha_A^*(\eta_1)\wedge d^0\overline{\varphi^*}(\eta_2)\Big).
\end{eqnarray*}
By (\ref{s34}),
\begin{eqnarray*}
&&d^0(\eta_1\wedge\eta_2)(X,Y,Z)\\
&=&\varphi^*\rho_A(X)\eta_1(\alpha_A(Y))\eta_2(\alpha_A(Z))-\varphi^*\rho_A(X)\eta_1(\alpha_A(Z))\eta_2(\alpha_A(Y))+\eta_1(\alpha_A(Y))\varphi^*\rho_A(X)\eta_2(\alpha_A(Z))\\
&&-\eta_1(\alpha_A(Z))\varphi^*\rho_A(X)\eta_2(\alpha_A(Y))-\varphi^*\rho_A(Y)\eta_1(\alpha_A(X))\eta_2(\alpha_A(Z))+\varphi^*\rho_A(Y)\eta_1(\alpha_A(Z))\eta_2(\alpha_A(X))\\
&&-\eta_1(\alpha_A(X))\varphi^*\rho_A(Y)\eta_2(\alpha_A(Z))+\eta_1(\alpha_A(Z))\varphi^*\rho_A(Y)\eta_2(\alpha_A(X))+\varphi^*\rho_A(Z)\eta_1(\alpha_A(X))\eta_2(\alpha_A(Y))\\
&&-\varphi^*\rho_A(Z)\eta_1(\alpha_A(Y))\eta_2(\alpha_A(X))+\eta_1(\alpha_A(X))\varphi^*\rho_A(Z)\eta_2(\alpha_A(Y))-\eta_1(\alpha_A(Y))\varphi^*\rho_A(Z)\eta_2(\alpha_A(X))\\
&&-\eta_1([X,Y])\eta_2(\alpha_A(Z))+\eta_1(\alpha_A(Z))\eta_2([X,Y])+\eta_1([X,Z])\eta_2(\alpha_A(Y))-\eta_1(\alpha_A(Y))\eta_2([X,Z])\\
&&-\eta_1([Y,Z])\eta_2(\alpha_A(X))+\eta_1(\alpha_A(X))\eta_2([Y,Z])\\
&=&\varphi^*\rho_A(X)\eta_1\wedge\eta_2(\alpha_A(Y),\alpha_A(Z))-\varphi^*\rho_A(Y)\eta_1\wedge\eta_2(\alpha_A(X),\alpha_A(Z))+\varphi^*\rho_A(Z)\eta_1\wedge\eta_2(\alpha_A(X),\alpha_A(Y))\\
&&-\eta_1\wedge\eta_2([X,Y],\alpha_A(Z))+\eta_1\wedge\eta_2([X,Z],\alpha_A(Y))-\eta_1\wedge\eta_2([Y,Z],\alpha_A(X)).
\end{eqnarray*}
Then, for any $\eta\in C^2(\Gamma(A);\CWM)$, we have
\begin{eqnarray*}
d^0\eta(X,Y,Z)&=&\varphi^*\rho_A(X)\eta(\alpha_A(Y),\alpha_A(Z))-\varphi^*\rho_A(Y)\eta(\alpha_A(X),\alpha_A(Z))+\varphi^*\rho_A(Z)\eta(\alpha_A(X),\alpha_A(Y))\\
&&-\eta([X,Y],\alpha_A(Z))+\eta([X,Z],\alpha_A(Y))-\eta([Y,Z],\alpha_A(X)).
\end{eqnarray*}
Sept5, for $\xi\in\Gamma(A^*)$, by(\ref{s34}), (\ref{s36}) and \rm{(i)}, we have:
\begin{eqnarray*}
0&=&d^0\circ d^0\xi(X,Y,Z)\\
&=&\varphi^*\rho_A(X)d^0\xi(\alpha_A(Y),\alpha_A(Z))-\varphi^*\rho_A(Y)d^0\xi(\alpha_A(X),\alpha_A(Z))+\varphi^*\rho_A(Z)d^0\xi(\alpha_A(X),\alpha_A(Y))\\
&&-d^0\xi([X,Y],\alpha_A(Z))+d^0\xi([X,Z],\alpha_A(Y))-d^0\xi([Y,Z],\alpha_A(X))\\
&=&\xi([[X,Y],\alpha_A(Z)]+[[Y,Z],\alpha_A(X)]+[[Z,X],\alpha_A(Y)]).
\end{eqnarray*}
So, we have:
\begin{equation}\label{s38}
[[X,Y],\alpha_A(Z)]+[[Y,Z],\alpha_A(X)]+[[Z,X],\alpha_A(Y)]=0.
\end{equation}

By (\ref{s35}) and (\ref{s38}), $(\Gamma(A),[\cdot,\cdot],\alpha_A)$ is a Hom-Lie algebra;

By(\ref{s37}), we have: $[X,fY]=\varphi^*(f)[X,Y]+\rho_A(X)f\alpha_A(Y)$;

By (\ref{s33}) and (\ref{s36}),  $\rho_A$ is a representation of Hom-Lie algebra $(\Gamma(
A),[\cdot,\cdot],\alpha_A)$ on $\CWM$ with respect to $\varphi^*$.

The proof is completed.\qed

\begin{thm}\label{thm5}
Let $B$ be a vector
bundle over  manifold $M, \varphi:M\longrightarrow M$ is a smooth
map and $\varphi^2=\rm{id}$, $\alpha_B:\Gamma(B) \longrightarrow
 \Gamma(B)$ is a linear endomorphism of $\Gamma(B)$ i.e.for $f\in\CWM,X\in\Gamma(B), \alpha_B(fX)=\varphi^*(f)\alpha_B(X)$. Then $(B,\varphi,[\cdot,\cdot],\rho_B,\alpha_B)$ is a Hom-Lie algebriod define by Definition\ref{s1-2} if and only if there is a
series operators $d^s:C^k(\Gamma(B);\CWM)\longrightarrow C^{k+1}(\Gamma(B);\CWM), s=0,1,\ldots$, and such that:
\begin{itemize}
\item[\rm{1)}]$d^s\circ d^s=0$;
\item[\rm{2)}]for any $\eta_1\in C^k(\Gamma(B);\CWM),\eta_2\in C^l(\Gamma(B);\CWM)$, we have
$$d^s(\eta_1\wedge\eta_2)=d^{s+l}\eta_1\wedge\overline{\varphi^*}\circ\alpha_B^*(\eta_2)+(-1)^k\overline{\varphi^*}\circ\alpha_B^*(\eta_1)\wedge
d^{s+k}\eta_2.$$
\item[\rm{3)}]$\alpha_B^*\circ d^s=d^s\circ\alpha_B^*;\quad \overline{\varphi^*}\circ d^s=d^{s+1}\circ\overline{\varphi^*}.$
\item[\rm{4)}]for $f\in\CWM=C^0(\Gamma(B);\CWM)$, we have: $d^1f\in \Gamma(B^*)$.
\item[\rm{5)}]for $f\in\CWM, \xi\in\Gamma(B^*)$, we have:
$$d^1\xi(fX_1,X_2)=\varphi^*(f)d^1\xi(X_1,X_2).$$
\end{itemize}
\end{thm}
\pf   For necessity, with above Propositions which we proved, we just need to prove \rm{4)} and \rm{5)}.

For Hom-Lie algebriod $(B,\varphi,[\cdot,\cdot],\rho_B,\alpha_B)$ and $f,g\in\CWM=C^0(\Gamma(B);\CWM), X\in\Gamma(B)$, by the definition of $d^s$, we have:
$$d^1f(gX)=\rho_B(gX)\varphi^*f=g\rho_B(X)\varphi^*f=gd^1f(X).$$

For Hom-Lie algebriod $(B,\varphi,[\cdot,\cdot],\rho_B,\alpha_B)$ and $f\in\CWM, \xi\in\Gamma(B^*)$, we have:
\begin{eqnarray*}
d^1\xi(fX_1,X_2)&=&\varphi^*\rho_B(fX_1)\xi(\alpha_B(X_2))-\varphi^*\rho_B(X_2)\varphi^*\xi(\alpha_B(fX_1))\\
&&-\xi([fX_1,X_2])\\
&=&\varphi^*(f)\varphi^*\rho_B(X_1)\xi(\alpha_B(X_2))-\varphi^*\Big(f\rho_B(X_2)\xi(\alpha_B(X_1))+\varphi^*\xi(\alpha_B(X_1))\rho_B(X_2)\varphi^*(f)\Big)\\
&&-\xi(\varphi^*(f)[X_1,X_2]-\rho_B(\alpha_B(X_2))f\alpha_B(X_1))\\
&=&\varphi^*(f)d^1\xi (X_1,X_2).
\end{eqnarray*}
So, we proved the necessity of this Theorem.  The sufficiency of this theorem is similar with Theorem \ref{thm4}.
\qed

\begin{rmk}
\begin{itemize}
\item[$\bullet$]
 For $E$ is a vector bundle over $M$, we hope that we are able to get a pull-back diagram:
\begin{equation*}
\xymatrix@C=0.5cm{
 E \ar[d]  \ar[rr] && \varphi^!E
                  \ar[d]  \ar[rr] && E  \ar[d]  \\
 M  \ar[rr]^{\varphi} &&
                M  \ar[rr]^{\varphi} && M  .
                }
\end{equation*}
So, we assume that $\varphi^2=\rm{id}$.
\item[$\bullet$]For given Hom-Lie algebriod $E$, we have operators $d^s$, then, by operators $d^s$, we can get the other definitions of Hom-Lie algebriod.
\end{itemize}
\end{rmk}

\end{document}